\documentclass[12pt]{article}

\usepackage{amsfonts,amssymb,amsmath,amsthm,eucal}
\numberwithin{equation}{section}

\newcommand{\la}{\lambda}
\newcommand{\sq}{{\square}}
\newcommand{\squ}{\sq_{\textup{\tiny up}}}
\newcommand{\sql}{\sq_{\textup{\tiny left}}}
\newcommand{\squl}{\sq_{\textup{\tiny up\&left}}}
\newcommand{\LL}{\mathsf{L}}
\newcommand{\Z}{\mathbb{Z}}
\newcommand{\Q}{\mathbb{Q}}

\DeclareMathOperator{\row}{row}
\DeclareMathOperator{\col}{col}
\DeclareMathOperator{\ord}{ord}
\DeclareMathOperator{\tr}{tr}
\DeclareMathOperator{\const}{const}

\newtheorem{Theorem}{Theorem}
\newtheorem{Proposition}{Proposition}
\newtheorem{Lemma}{Lemma}

\begin{document}

\title{Combinatorial formula for Macdonald polynomials, Bethe Ansatz,
and generic Macdonald polynomials}
\author{Andrei Okounkov\thanks{
 Department of Mathematics, University of California at
Berkeley, Evans Hall \#3840, 
Berkeley, CA 94720-3840. E-mail: okounkov@math.berkeley.edu}
}
\date{}

\maketitle

\begin{abstract} 
We give a direct proof of the combinatorial formula
for interpolation Macdonald polynomials by introducing
certain polynomials, which we call generic Macdonald 
polynomials, which depend on $d$ additional parameters and specialize
to all Macdonald polynomials of degree $d$. The form of these
generic polynomials is that of a Bethe eigenfunction and they
imitate, on a more elementary level,
the $R$-matrix construction of quantum immanants. 
\end{abstract}

\section{Introduction}

\subsection{Interpolation Macdonald polynomials and combinatorial formula}

\subsubsection{The polynomials $I_\mu$} 

Consider the symmetric polynomials
$$
I_\mu(x_1,\dots,x_n;q,t,s)\in \Q(q,t,s)[x_1,\dots,x_n]^{S(n)}\,,
$$
where $\mu$ is a partition with at most $n$ parts, satisfying
the following interpolation condition: for any partition 
$\la$ such that $\mu\not\subset \la$
\begin{equation}\label{van}
I_\mu(\boldsymbol{\la})=0\,, \quad
\boldsymbol{\la}=([1,\la_1],\dots,[n,\la_n])\,,
\end{equation}
where
$$
[i,j]=\frac{sq^{j+1}}{t^{i+1}}+\frac{t^{i+1}}{sq^{j+1}}\,. 
$$
Together with the condition that the degree of $I_\mu$
is at most $|\mu|$ this determines $I_\mu$ uniquely up
to a constant factor. These polynomials $I_\mu$ were
introduced and studied in \cite{Ok4}; they
specialize to
interpolation polynomials studied by Knop, Olshanski,
Sahi, and the author in a series of papers, see
the References. We will call the polynomials $I_\mu$
the interpolation Macdonald polynomials. 

\subsubsection{Combinatorial formula} 
The polynomials $I_\mu$ are very distinguished
multivariate special functions with a large number of
deep properties and important applications (in particular,
 to the better known orthogonal Macdonald polynomials). Most of 
their properties can be quite easily deduced from the following 
explicit formula for these polynomials which is known as the
combinatorial formula. 

Consider the diagram of the partition $\mu$ and a 
function 
$$
f:\sq \mapsto \{1,\dots,n\}
$$
from the squares of this diagram to integers. The function $f$
is called a \emph{reverse tableau} if its values weakly
decrease along the rows and strictly decrease along the 
columns of $\mu$. We have the following formula \cite{Ok4}
\begin{equation}\label{combfor0}
I_\mu=\sum_{\textup{reverse tableux $f$}} \Psi_f \prod_{\sq\in\mu}
\left(x_{f(\sq)} - \frac{s q^{\col(\sq)}}{t^{f(\sq)+\row(\sq)}}
-
\frac{t^{f(\sq)+\row(\sq)}}{s q^{\col(\sq)}}
\right) \,,
\end{equation}
where $\col(\sq)$ and $\row(\sq)$ denote the the column and
row numbers a square $\sq\in\mu$ and $\Psi_f$ is a certain
weight of the tableau $f$ which will be specified below.
It is the same weight as appears in the combinatorial
formula for the ordinary Macdonald polynomials, which makes
it clear that the ordinary Macdonald polynomial is the top
degree term of $I_\mu$. 

\subsubsection{Direct proof of the combinatorial formula}

One of the main results of this note is a proof of 
the combinatorial formula 
\eqref{combfor0} which does not use any properties
of $I_\mu$ except for their definition. We recall that a
great deal of nontrivial properties of $I_\mu$
were used in \cite{Ok4} in the proof of \eqref{combfor0}.
Now many of  these properties can now be without difficulty deduced 
from  \eqref{combfor0} without falling into a vicious 
circle. 

It is obvious that the right-hand side of \eqref{combfor0} is a polynomial
of degree $|\mu|$ and it is also easy to check that
it satisfies the vanishing condition \eqref{van},
see for example Lemma 4.1 in \cite{Ok5}. What remains is
the nontrivial task of checking that this polynomial is
symmetric in $x_1,\dots,x_n$. 

\subsubsection{How we prove the symmetry}

Our way of verifying the symmetry in \eqref{combfor0}
will be the following. We introduce certain polynomials
$$
I_U(x)\in\Q(u_1,\dots,u_d,q,t,s)[x_1,\dots,x_n]
$$
of degree $d$ in $x$. This polynomials have the property that
for any $\mu$ such that $|\mu|=d$ one can tune the parameters
$\{u_i\}$ so that to make $I_U$ go into
the right-hand side of \eqref{combfor0}. Because of this property,
we call the polynomials $I_U$ the \emph{generic
Macdonald polynomials}.

The additional freedom which comes with the parameters
$\{u_i\}$ makes it easy to check that the polynomials
$I_U(x)$ are symmetric. This is done in Section \ref{s2}.
After that, in Section \ref{s3} we show how to specialize
the polynomials $I_U$ to polynomials in the right-hand
side of \eqref{combfor0} which, in particular, completes
our new proof of \eqref{combfor0}. 

\subsection{Generic Macdonald polynomials}

\subsubsection{Shifted Schur functions and Quantum Immanants}

The polynomials $I_U$ are not only a convenient
technical tool but also an interesting object of study 
in their own right.

In order to better motivate their definition,  
we begin with the special case when $q=t\to 1$ and 
$s\to\infty$. In this limit, the polynomials $I_\mu$ become (after a
change of variables)  the shifted
Schur functions $s^*_\mu$, see \cite{OO}, for which the 
combinatorial formula is 
\begin{equation}\label{cfs}
s^*_\mu(x)=\sum_{\textup{reverse tableux $f$}} \,\, \prod_{\sq\in\mu}
\left(x_{f(\sq)} - c(\sq)\right)\,,
\end{equation}
where 
$$
c(\sq) = \col(\sq)  - \row(\sq) 
$$ 
denotes the \emph{content} of a square $\sq\in\mu$ \,.

The functions $s^*_\mu$  have a very nice and important
representation-theoretic 
interpretations as the image under the Harish-Chandra
homomorphism of a distinguished linear basis of the
center of the universal
enveloping algebra $\mathcal{U}(\mathfrak{gl}_n)$.
The corresponding central elements are known
as the \emph{quantum immanants} and denoted 
by $\mathbb{S}_\mu$.  

Explicit formulas for quantum immanants were obtained in
\cite{Ok1}, see also \cite{N,Ok11}. 
These formulas can be viewed as a noncommutative analog,
or quantization, of the combinatorial formula \eqref{cfs},
see for example Section 3.7 of \cite{Ok1}.
We recall from \cite{N} the following
construction of $\mathbb{S}_\mu$ which was motivated by
representation theory of the Yangian $Y(\mathfrak{gl}_n)$.

Consider the algebra
of matrices with entries in  $\mathcal{U}(\mathfrak{gl}_n)$.
Let $E=(E_{ij})$ be the matrix formed by
the standard basis elements of $\mathfrak{gl}_n$. Denote
$|\mu|$ by $d$ and consider the following expression
\begin{equation}\label{parovoz}
\tr\, (E-u_1)\otimes(E-u_2)\otimes \cdots \otimes(E-u_d) \cdot
\overrightarrow{\prod_{i,j}}\, \mathsf{R}_{ij}(u_i-u_j)  \,,
\end{equation}
where 
$$
\mathsf{R}_{ij}(u) = 1 - \frac{(ij)}{u} \in \Q(u)\, S(d)
$$
is the rational $R$-matrix, the ordered product in \eqref{parovoz}
ranges 
over all pairs $1\le i<j\le d$ ordered lexicographically, the 
symmetric group $S(d)$ acts on the $d$-fold tensor product
by the usual permutation of indices,  and the 
trace means, as usual, the sum of the diagonal matrix elements.

The expression \eqref{parovoz} is an element of
$$
\mathcal{U}(\mathfrak{gl}_n) \otimes \Q(u_1,\dots,u_d)  \,.
$$
It is explained in \cite{N} how \eqref{parovoz} 
turns into $\mathbb{S}_\mu$ when the parameters
$\{u_i\}$ approach, in a certain special way, 
the contents of the diagram $\mu$. 
In this limit, the product of $R$-matrices in \eqref{parovoz}  
turns into a projection onto a certain specific vector in the  irreducible 
representation of $S(d)$ corresponding to $\mu$.

\subsubsection{Formula for generic Macdonald polynomials}

The definition of the polynomials $I_U$ is our attempt to
keep as much as possible from \eqref{parovoz} in the
case of general $q$ and $t$. This general case has two
main features. On the one hand, in the absence
of a simple structure to replace $\mathcal{U}(\mathfrak{gl}_n)$ for the 
general $q$ and $t$, we do not go over the Harish-Chandra bridge and
stay entirely within the commutative world
of polynomials. On the other hand, the additional
flexibility provided by parameters $q$ and $t$ allows us
to work with polynomial in a way very similar to
what was done in the noncommutative world of $\mathcal{U}(\mathfrak{gl}_n)$. 

As a commutative replacement for the $R$-matrix $\mathsf{R}(u)$, we introduce
the following function 
$$
\rho(u)=\frac{(1-qu)(1-tu/q)}{(1-u)(1-tu)} 
$$
and define the polynomial $I_U$ to be the following analog of
\eqref{parovoz}
\begin{equation}\label{Rsum}
I_U(x_1,\dots,x_n)=\sum_{i_1,\dots,i_d=1}^n \, 
\prod_{k=1}^d 
\left(x_{i_k} + \frac{u_k}{t^{i_k}} +  \frac{t^{i_k}}{u_k}
\right) \prod_{i_k< i_l }
\rho(u_k/u_l)\,,
\end{equation}
where the summation replaces the trace,
first product replaces the tensor product, and the second
product replaces the product of $R$-matrices. 

Of course, the general structure of the formula \eqref{Rsum} is
a typical form of a Bethe  eigenfunction.

\subsubsection{Properties of generic Macdonald polynomials}

As already mentioned, in Theorem \ref{t1} below we prove that 
polynomials $I_U$ are symmetric in $x_1,\dots,x_n$. The
proof of the symmetry is elementary and based on the
study of the singularities of $I_U$ in parameters $\{u_i\}$.

In Section \ref{s3} we show that  as the parameters $\{u_i\}$
approach the $(q,t)$-analogs of the contents of a diagram $\mu$, that is, as
\begin{equation}\label{predel}
\{u_i\} \to \left\{\const q^{\col(\sq)} \, t^{-\row(\sq)}\right\}_{\sq\in\mu}
\end{equation}
the generic Macdonald polynomial becomes the interpolation
Macdonald polynomial labeled by $\mu$, see Theorem \ref{t2}. 
The limit \eqref{predel} has to be taken
in a certain special way which is parallel to what one
does to obtain $\mathbb{S}_\mu$ as a limit of
\eqref{parovoz}. 

Finally, we mention an important but potentially difficult open
problem which is to 
expand the generic polynomial $I_U$ in terms of the
polynomials $I_\mu$.

\section{The polynomials $I_U(X)$}\label{s2} 

\subsection{Definition}\label{def}
Let $X=\{x_1,\dots,x_n\}$ and $U=\{u_1,\dots,u_d\}$ be two finite
sets of variables. Define a polynomial $I_U(x)$ of degree $d=|U|$ 
in the variables $X$ with coefficients in rational functions in $U$ as follows. 

Introduce the following function
$$
\rho(v)=\frac{(1-qv)(1-tv/q)}{(1-v)(1-tv)} \,,
$$
where $q$ and $t$ are parameters. Given a function $f$  
$$
f: \{1,\dots,d\} \to \{1,\dots,n\}
$$
set, by definition,
$$
R_f= \prod_{f(i)<f(j)} \rho(u_i/u_j) \,,
$$  
where the product is over all pairs $i,j$ in $\{1,\dots,d\}$ such that
$f(i)<f(j)$. The polynomial $I_U(x)$ is defined by the following
formula
\begin{equation}\label{definition}
I_U(X)=\sum_{f} R_f\, \prod_{i=1}^d 
\left(x_{f(i)} + \frac{u_i}{t^{f(i)}} + \frac{t^{f(i)}}{u_i}\right)  \,,
\end{equation}
where the summation is over all $n^d$ possible functions $f$. It is clear
that the definition \eqref{definition} is just another way of
writing the formula \eqref{Rsum} from the Introduction. We call
these polynomials $I_U$ the \emph{generic Macdonald polynomials}.

{}From definition of $I_U(X)$, the following 
decomposition property of this
polynomial is obvious. Fix any $k=1,\dots,n$ and set
$X'=\{x_1,\dots,x_k\}$. 
Then
\begin{equation}\label{decomposition}
I_U(X)=\sum_{U'\subset U} I_{U'}(X')\,  I_{U\setminus U'} (X\setminus X') \!\!
\prod_{u_i\in U', u_j\in U\setminus U'} \!\!\rho(u_i/u_j) \,,
\end{equation}
where the summation is over all subsets $U'\subset U$. 

\subsection{Symmetry}

Since no ordering on $U$ was used, the polynomial $I_U(X)$ is 
obviously invariant under permutations of the $u_i$'s. 
Our main result about the polynomials $I_U(X)$ is that it is
also invariant under permutations of the $x_i$'s. 

\begin{Theorem}\label{t1}
The polynomial $I_U(X)$ is symmetric in $x_1,\dots,x_n$.
\end{Theorem} 

In the proof, we shall need some elementary properties of the
function $\rho(v)$ which we collect in the following 

\begin{Lemma} We have
\begin{align}
\rho(v)&=\rho\left(\frac1{tv}\right)\,, \label{recip}\\
\rho(v)&=1+\frac{(q-1)(q-t)}{q} v + O\left(v^2\right)\,, \quad v\to 0\,,
\label{ser0}\\
\rho(v)&=1+\frac{(q-1)(q-t)}{qt} \frac1v + O\left(\frac1{v^2}\right)\,, \quad v\to \infty\,.
\label{serinf}
\end{align}
\end{Lemma}

\begin{proof}[Proof of Theorem \ref{t1}] 
It suffices to check that $I_U(X)$ is invariant under a
transposition of two adjacent $x_i$'s. By the
decomposition property \eqref{decomposition}, this reduces to showing
that $I_U(x_1,x_2)$ is symmetric in $x_1$ and $x_2$,  
identically in $u_1,\dots,u_d$ for any $d$. 

Let us single out
one of the $u_i$'s  and write
$$
g(u_1)= I_U(x_1,x_2)-I_U(x_2,x_1)\,,
$$
where $x_1$, $x_2$, and the rest of the $u_i$'s are considered as
parameters. Our goal is to show that $g=0$ identically.

The function $g$ is a rational function of $u_1$ with at most
first order poles at the following points 
$$
u_1 \in \{0,\infty,u_i, t^{\pm1} u_i\}\,, \quad i=2,\dots,d\,.
$$ 
First we show that, in fact, this function is regular everywhere. 
By symmetry,
it suffices to consider the points $\{0,\infty,u_2,t^{\pm1} u_2\}$. 

Observe that on the divisor $u_1=u_2$ the function $I_U(x_1,x_2)$ 
itself is regular. Indeed, a function which is symmetric in $u_1$ and
$u_2$ cannot have a pole of exact order one on the diagonal $u_1=u_2$. 

Next, consider the divisor $u_1=t^{-1} u_2$. Only those functions $f$
which satisfy
$$
f(1)=1\,, \quad f(2)=2
$$ 
contribute to the pole of $I_U(x_1,x_2)$ 
at $u_1=t^{-1} u_2$. For such function $f$ and 
any $i=3,\dots,d$ we have either 
$f(i)=f(1)<f(2)$ or $f(1)<f(2)=f(i)$. 
The condition $u_1=t^{-1} u_2$ and \eqref{recip} yield that 
$$
\rho\left(\frac{u_i}{u_2}\right)=\rho\left(\frac{u_1}{u_i}\right) \,.
$$  
Therefore
$$
R_f = R_{f|_{\{3,\dots,d\}}} \, \prod_{i=2}^d \rho\left(\frac{u_1}{u_i}\right) \,,
$$
where $f|_{\{3,\dots,d\}}$ denotes the restriction of $f$ to $\{3,\dots,d\}$. 
It follows that 
the residue of $I_U(x_1,x_2)$ at $u_1=t^{-1} u_2$ is proportional to
$$
\left(x_1+\frac{u_1}{t}+\frac{t}{u_1}\right)
\left(x_2+\frac{u_2}{t^2}+\frac{t^2}{u_2}\right)
 \, I_{\{u_3,\dots,u_d\}}(x_1,x_2)\,,
$$
which by the condition $u_1=t^{-1} u_2$ and induction on the number of $u_i$'s is 
symmetric in $x_1$ and $x_2$. 

In other words, $g(u_1)$ is regular at $u_1=t^{-1} u_2$. By the symmetry
between $u_1$ and $u_2$, it is also regular at $u_1=t u_2$.

Now consider the point $u_1=0$. 
Since for any function $f$ we either have $f(1)=1$ or $f(1)=2$,
it follows from the definition of $I_U(x_1,x_2)$ that
\begin{multline}
I_U(x_1,x_2) =\sum_{f'} R_{f'} \Pi_{f'} \left[  \left(x_1+\frac{u_1}{t}+\frac{t}{u_1}\right)
\prod_{f'(i)=2} \rho  \left(\frac{u_1}{u_i}\right)+ \right. \\
\left. \left(x_2+\frac{u_1}{t^2}+\frac{t^2}{u_1}\right) 
\prod_{f'(i)=1} \rho  \left(\frac{u_i}{u_1}\right)\right]\,,
\end{multline}  
where the summation is over all possible function 
$$
f':\{2,\dots,d\} \to \{1,2\} \,,
$$
and where we use the following abbreviation
$$
\Pi_{f'} = \prod_{i=2}^d \left(x_{f'(i)} + \frac{u_i}{t^{f'(i)}} + \frac{t^{f'(i)}}{u_i}\right)\,.
$$
{}From \eqref{ser0} we have
\begin{multline*}
\left(x_1+\frac{u_1}{t}+\frac{t}{u_1}\right) \prod_{f'(i)=2} \rho  
\left(\frac{u_1}{u_i}\right) = \\
\frac{t}{u_1} + x_1+ \frac{t(q-1)(q-t)}{q} \sum_{f'(i)=2} \frac{1}{u_i} + 
O(u_1) \,, \quad u_1\to 0 \,.
\end{multline*}
Similarly, from \eqref{serinf} we obtain 
\begin{multline*}
\left(x_2+\frac{u_1}{t^2}+\frac{t^2}{u_1}\right) \prod_{f'(i)=1} \rho  
\left(\frac{u_i}{u_1}\right) =\\
\frac{t^2}{u_1} + x_2 + \frac{t(q-1)(q-t)}{q} \sum_{f'(i)=1} \frac{1}{u_i}
 + O(u_1) \,, \quad u_1\to 0 \,.
\end{multline*}
Therefore,
\begin{multline*}
I_U(x_1,x_2)= \\
\left(\frac{t+t^2}{u_1}+x_1+x_2+ \frac{t(q-1)(q-t)}{q}\sum_{i=2}^d \frac{1}{u_i} \right)
I_{\{u_2,\dots,u_d\}}(x_1,x_2) + O(u_1)\,,
\end{multline*}
as $u_1\to 0$, which means that both the $u_1^{-1}$ and $u_1^0$ terms are symmetric in
$x_1$ and $x_2$. In other words, not only is $g(u_1)$ regular at $u_1=0$ but it also
vanishes at $u_1=0$. Similarly, $g(u_1)$ is regular at $\infty$ and thus, as a function
 regular everywhere and vanishing at one point, $g$ is zero identically. 
\end{proof}

\section{Specialization to interpolation Macdonald polynomials}\label{s3}

\subsection{Limit transition}
Let $\la$ be the diagram of a partition of $d$. The diagram $\lambda$
has $d$ squares and it will be convenient to consider the set $U$ to
be indexed by the squares of $\la$ rather than integers $\{1,\dots,d\}$. 
Denote by $I_\la(X;q,t,s,r)$ the image of
$I_U(X)$ under the following
specialization of the variables $U$
\begin{equation}\label{spec}
u_\sq = - s\, q^{\col(\sq)} r^{\row(\sq)}\,, \quad \sq\in\la \,.
\end{equation}
Here $s$ and $r$ are parameters,
$\sq\in\la$ is a square in the diagram $\la$, and 
$\col(\sq)$ and $\row(\sq)$ are the column and row numbers of $\sq$,
respectively. Also, 
the functions $f$ will now be functions from the diagram $\la$ to $\{1,\dots,n\}$.
Such functions are known in combinatorics as \emph{tableaux} on $\la$. 

Let $R_f(\la;q,t,r)$ be the image of $R_f$ under the specialization \eqref{spec}.
Because only the ratios of the $u$'s enter the formula for $R_f$, this specialization
does not depend on $s$. 
This image is nontrivial only for some special $f$'s as the following proposition
shows: 

\begin{Proposition}\label{p1} We have
$$
R_f(\la;q,t,r)=0
$$
unless $f$ is row decreasing, that is, unless the values of $f$ weakly decrease
inside the rows. 
\end{Proposition}

\begin{proof} Suppose $f$ is not row decreasing. Then there 
exist two
squares $\sq$ and $\sql$ in $\la$ such that $\sql$ is just left of $\sq$
and
$$
f(\sql) < f(\sq) \,.
$$
In this case, $R_f(\la;q,t,r)$ has the factor of
$$
\rho\left(\frac{u_{\sql}}{u_{\sq}}\right) =
\rho\left(\frac1q\right)=0 \,,
$$
and hence $R_f(\la;q,t,r)=0$.
\end{proof}

The main result of this Section is the following:

\begin{Theorem}\label{t2} We have
\begin{equation}\label{lim}
I_\la(X;q,t,s) \propto \lim_{r\to t^{-1}} (1-r t)^{\ell(\la)-1} I_\la(X;q,t,s,r) \,,
\end{equation}
where $I_\la(X;q,t,s)$ is the interpolation Macdonald polynomial
corresponding to the partition $\la$ and the proportionality constant
lies in $\Q(q,t)$. 
\end{Theorem}

Theorem \ref{t2} will be established by comparing the limit
in \eqref{lim} with the right-hand side of \eqref{combfor0},
which in view of the symmetry of the $I_U$ established above,
simultaneously proves \eqref{lim} and the combinatorial formula. 

First, we check that only those $f$ which are reverse tableaux
contribute to the limit \eqref{lim}. Then, in the next subsection, we shall
check that for such functions $f$ the weights 
\begin{equation}\label{e1}
R_f(\la;q,t)= \lim_{r\to t^{-1}} (1-r t)^{\ell(\la)-1} R_f(\la;q,t,r)
\end{equation}
are proportional to the weights $\Psi_f$
in the combinatorial formula \eqref{combfor0}. 

\begin{Proposition}\label{p2} We have 
$$
R_f(\la;q,t)=0
$$
unless $f$ is a reverse tableaux on $\la$. If $f$ is a reverse tableaux
then $R_f(\la;q,t)$ is finite and nonzero. 
\end{Proposition} 

\begin{proof}[Proof of Proposition \ref{p2}]
Recall that $\rho(v)$ has two poles $v=1,t^{-1}$ and two
zeros $v=q^{-1},qt^{-1}$. The pole $v=1$ is not an issue
in limit \eqref{e1} and the zero $v=q^{-1}$ has been
already accounted for in Proposition \ref{p1}. The pole
$v=t^{-1}$ means that
$$
\rho\left(\frac{u_{\sq}}{u_{\squ}}\right) \to \infty\,, \quad r\to\frac 1t\,,
$$
where  $\squ$
is the square right on top of the square $\sq$. This factor is present
in  $R_f(\la;q,t,r)$ if 
$$
f(\sq) < f(\squ)\,.
$$
Similarly, the zero
$v=qt^{-1}$ means that
$$
\rho\left(\frac{u_{\sq}}{u_{\squl}}\right) \to 0\,, \quad r\to\frac 1t\,,
$$
where $u_{\squl}$ is the square to the left of $\squ$. This
factor enters $R_f(\la;q,t,r)$ if 
$$
f(\sq) < f(\squl)\,.
$$

Let $\ord_f$ be the order of the pole of $R_f(\la;q,t,r)$ at
$r=t^{-1}$. We will show that
\begin{equation}\label{e2}
\ord_f \le \ell(\la)-1
\end{equation}
and that this maximum is reached
precisely for reverse tableaux. We have 
\begin{equation}\label{ord}
\ord_f = \#\{\sq,f(\sq) < f(\squ)\} - \#\{\sq,f(\sq) < f(\squl)\}  \,.
\end{equation}
The contribution of any square $\sq\in\la$ to \eqref{ord} is $0$ or $\pm1$.

If the square $\sq$ is not in the first column of $\la$
then its contribution to \eqref{ord} is $\le 0$. Indeed,
if $f(\sq) < f(\squ)$ then since $f$ is row decreasing we
also have 
$$
f(\sq) < f(\squ) \le f(\squl)\,.
$$
That is, only the squares in the first column can make positive 
contribution to \eqref{ord} which proves \eqref{e2}.
It is also clear that if $f$ is a reverse tableaux then
$\ord_f=\ell(\la)-1$. 

Now suppose that $f$ is not a reverse tableaux. In order to
prove that
$$
\ord_f<\ell(\la)-1
$$
it suffices to find one square whose contribution to $\ord_f$
is negative. By assumption,
there exists a square $\sq$ such that $f(\sq) \ge  f(\squ)$.
Additionally, we can assume that such a $\sq$ is minimal in the sense that
$f(\sql) <  f(\squl)$. Then since $f$ is row decreasing we
conclude that $f(\sq) \le f(\sql) <  f(\squl)$ which means
that the contribution
of this square $\sq$ to \eqref{ord} equals $-1$. It follows that
$\ord_f<\ell(\la)-1$. This concludes the proof. 
\end{proof} 

\subsection{Computation of $R_f(\la;q,t)$}

\subsubsection{The map $\LL$} 

The weights $R_f(\la;q,t)$ are products of a large number
of factors of the form $(1-q^k t^l)$ where $k,l\in\Z$. We
will treat such factors as single indivisible objects, ignoring
the fact that $(1-q^k t^l)$ is reducible if $\gcd(k,l)>1$.  

Therefore, it will be convenient to switch, following Macdonald, 
from the multiplicative to the additive
notation as follows. Introduce a map $\LL$ from the free Abelian group
(written multiplicatively) generated by the symbols $(1-a)$ to the
free Abelian group (written additively) generated by symbols $a$ 
which is defined on generators by 
$$
\LL((1-a))=a \,.
$$

The map $\LL$ is also a very convenient way to deal with zeros and poles.
If in an expression $A(z)$ depending on a parameter $z$ 
$$
A(z)=\prod_i (1-a_i(z))^{\pm 1}
$$
$m_+$ parentheses in the numerator and $m_-$ parenthesis in the
denominator vanish at a certain point $z=z_0$
then $\LL(A(z))$ will have a summand of $m_+-m_-$ at $z=z_0$.   

We agree that
$$
\LL^{-1} \left(\frac{x}{1-q}\right) = \prod_{i=0}^\infty (1-q^i x) \,,
$$
which converges as a formal series in $q$ or else converges in the usual analytic
sense if $|q|<1$. Normally, we will only encounter such linear
combinations of the fractions $\frac{x}{1-q}$ which are, in fact,
polynomials. 

Taking into consideration the relation $(1-a)=-a(1-a^{-1})$ we see
that it is permissible to  make the following transformation
\begin{equation}\label{transf}
\sum\pm a_i \mapsto \sum \pm a_i^{-1}
\end{equation}
provided $\prod (-a_i)^{\pm 1}= 1$,
where the choice of signs is the same as in \eqref{transf}.

Finally, observe that
\begin{align}
\LL(\rho(x)) = x\, \frac{(q-1)(q-t)}q \,, \label{e3} \\
\LL(\rho(x^{-1})) = x\, \frac{(q-1)(q-t)}{qt} \,. \label{e4}
\end{align}
These two formulas are related by the transformation \eqref{transf}. 

\subsubsection{Combinatorial formula for interpolation Macdonald
polynomials}

The combinatorial formula for the interpolation Macdonald polynomials \cite{Ok4}
expresses the polynomial corresponding to the partition $\la$ as a 
sum over all tableaux $f$ on $\la$ as follows
\begin{equation}\label{combfor}
I_\mu=\sum_{\textup{reverse tableux $f$}} \Psi_f \prod_{\sq\in\la}
\left(x_{f(\sq)} - \frac{s q^{\col(\sq)}}{t^{f(\sq)+\row(\sq)}}
-
\frac{t^{f(\sq)+\row(\sq)}}{s q^{\col(\sq)}}
\right) \,,
\end{equation}
Here $\Psi_f$ are a certain weight of associated to $f$ which will
be defined momentarily. 
Note that \eqref{combfor}
differs by a change of variables from the formula established
in \cite{Ok4}.

Recall that a reverse tableaux on $\la$ with entries in $\{1,\dots,n\}$
can be encoded as a sequence
of diagrams
\begin{equation}\label{chain}
\la=\la_0 \supset \la_1 \supset \la_2 \supset \dots \supset \la_n=\emptyset \,,
\end{equation}
such that each skew shape $\la_k/\la_{k+1}$ is a horizontal strip. This
encoding is simply the following: 
$$
\la_k=f^{-1}(\{k+1,\dots,n\}) \,.
$$
By definition,
$$
\Psi_f = \prod_{k=0}^{n-1} \psi_{\la_k/\la_{k+1}} \,,
$$
where (see Example 2(b) in Section VI.7 of \cite{M})
$$ 
\LL(\psi_{\la/\mu})=\frac{t-q}{1-q}\sum_{1\le i \le j \le \ell(\mu)}
\left(
q^{\mu_i-\mu_j}-q^{\la_i-\mu_j}-q^{\mu_i-\la_{j+1}}+q^{\la_i-\la_{j+1}}
\right) t^{j-i} \,.
$$
One can check that this is in fact a polynomial in $q$ and $t$.

\subsection{The weights $R_f(\la;q,t)$}

It is clear from their construction
that the weights $R_f(\la;q,t)$ factor similarly to the factorization
of $\psi_{\la/\mu}$
$$
R_f(\la;q,t)= \prod_{k=0}^{n-1} \rho_{\la_k/\la_{k+1}} \,,
$$
where
\begin{multline}\label{zetalamu}
\LL(\rho_{\la/\mu}) =\\
\frac{(q-1)(q-t)}q  
\sum_{\sq\in\mu,\,\boxtimes\in\la/\mu}  
q^{\col(\boxtimes)-\col(\sq)}\, t^{\row(\sq)-\row(\boxtimes)} +\ell(\la)-\ell(\mu)\,.
\end{multline}
Here the number $\ell(\la)-\ell(\mu)\in\{0,1\}$ is precisely the compensation for 
the singularity of $R_f(\la;q,t,r)$ as $r\to t^{-1}$. 

\begin{Proposition} The ratio $R_f(\la;q,t)\big/\Psi_f$ depends only on the
diagram $\la$ and not on the particular choice of a tableau $f$ on $\la$.
\end{Proposition}

It is clear that this proposition is all what remains to prove in order
to complete proof of Theorem \ref{t2}.

\begin{proof}

Let us split the sum in \eqref{zetalamu} in two 
$$
\LL(\rho_{\la/\mu})=\Sigma_+ + \Sigma_-\,,
$$
where $\Sigma_+$
and $\Sigma_{-}$ contains summands with nonnegative and negative number 
$\row(\sq)-\row(\boxtimes)$, respectively. 

Using the identity
$$
\sum_{k=a}^b \sum_{l=c}^d q^{k-l} = - \frac{q}{(1-q)^2} 
\left( q^{a-c} + q^{b-d} - q^{a-d-1} - q^{b-c+1} 
\right)
$$
to perform the summation along the rows, we compute:
$$
\Sigma_+ = \frac{q-t}{1-q} \sum_{i\le j \le \ell(\mu)} 
\left( q^{\la_i-\mu_j} + q^{\mu_i} - q^{\la_i} -q^{\mu_i-\mu_j} \right) t^{j-i} \,.
$$ 
Similarly, we compute $\Sigma_-$ using the formula \eqref{e4} instead of
\eqref{e3} and we find that
$$
\Sigma_- =  \frac{q-t}{1-q}  \sum_{i\le j \le \ell(\mu)} 
\left( q^{\mu_i-\la_{j+1}} + q^{-\mu_{j+1}} - q^{-\la_{j+1}} -q^{\mu_i-\mu_{j+1}} \right) t^{j-i} \,.
$$ 
For any diagram $\eta$ define $\xi_\eta$ by 
$$
\LL(\xi_{\eta}) =  \frac{q-t}{1-q} \sum_{i\le j\le \ell(\eta)} 
\left( q^{\eta_i} + q^{-\eta_{j+1}} - q^{\eta_i - \eta_{j+1}} \right)
t^{j-i} - \ell(\eta)+1 \,.
$$
Here the constant $1-\ell(\eta)$ compensates for the appearances of 
$q^0 t^0$ in the expansion of 
$\frac{q-t}{1-q} q^{-\eta_{j+1}} t^{j-i}$ which occur
precisely when $i=j < \ell(\eta)$. Also set by definition 
$$
\tau_{\la/\mu} =  \frac{\rho_{\la/\mu}}{\psi_{\la/\mu}}\frac{\xi_{\la}}{\xi_{\mu}} \,.
$$
It follows trivially from this definition that
\begin{equation}\label{tau}
\frac{R_f(\la;q,t)}{\Psi_f} = \frac{\xi_{\emptyset}}{\xi_{\lambda}} 
\prod_{k=0}^{n-1} \tau_{\la_k/\la_{k+1}}
\end{equation}

It is clear that 
$$
\LL(\tau_{\la/\mu}) =
   \frac{q-t}{1-q} \sum_{\ell(\mu)<j\le \ell(\la)} \sum_{i\le j} 
\left( q^{\la_i} + q^{-\la_{j+1}} - q^{\la_i - \la_{j+1}} \right) t^{j-i} \,. 
$$
Note that by our construction we have $\ell(\la)\le \ell(\mu)+1$. Therefore the above
sum is either zero or else the only possible value of $j$ in it is $j=\ell(\la)$.
In the latter case we have $\la_{j+1}=0$ which means that 
$$
\LL(\tau_{\la/\mu}) = \frac{q-t}{1-q} \sum_{i<\ell(\la)} t^{\ell(\la)-i} = 
\dfrac{q-t}{1-q} \dfrac{t^{\ell(\la)}-1}{t-1} \,.
$$
In other words, 
$$
\LL(\tau_{\la/\mu}) = \begin{cases}
0\,, & \ell(\la)=\ell(\mu)\,, \\ 
\dfrac{t-q}{1-q} \dfrac{t^{\ell(\la)}-1}{t-1}\,,   & \ell(\la)=\ell(\mu)+1\,.
\end{cases}
$$
It now clear that
$$
\prod_{k=0}^{n-1} \tau_{\la_k/\la_{k+1}}
$$
depends only on the first and last element in the chain \eqref{chain} which precisely
means that \eqref{tau} does not depend on the particular choice of the tableau $f$ 
on $\la$. 
\end{proof}


\begin{thebibliography}{99}





\bibitem{KOO}
 S.~Kerov, A.~Okounkov, and G.~Olshanski,
\emph{The boundary of Young graph with Jack edge multiplicities},
Internat.\ Math.\ Res.\ Notices 1998, no.\ 4, 173--199.




\bibitem{Kn}
 F.~Knop, \emph{ Symmetric and non--symmetric quantum Capelli
polynomials},
Comment.\ Math.\ Helv.\ 
\textbf{72}, 1997, 84--100.




\bibitem{KS}
 F.~Knop and S.~Sahi,
\emph{ Difference equations and symmetric polynomials
defined by their zeros}, 
 Internat.\ Math.\ Res.\ Notices, 
1996, no.\ 10, 473--486.






\bibitem{M}
 I.~G.~Macdonald,
\emph{Symmetric functions and Hall polynomials},
2nd ed.,
Oxford University Press, 1995.


\bibitem{N}
M.~L.~Nazarov,
\emph{Yangians and Capelli identities},
A.~A.~Kirillov Seminar on Representation Theory, 
edited by  G.~Olshanski,  
American Mathematical Society Translations, Ser.~2,
Amer.\ Math.\ Soc., Providence, 1997. 




\bibitem{Ok1}
 A.~Okounkov,
\emph{Quantum immanants and higher Capelli identities}
 Transformation Groups 
\textbf{1}, no.\ 1-2, 1996 , 99--126. 


\bibitem{Ok11}
 A.~Okounkov,
\emph{Young basis, Wick formula, and higher Capelli identities},
Internat.\ Math.\ Res.\ Notices 1996, no.\ 17, 817--839.




\bibitem{Ok2}
A.~Okounkov,
\emph{(Shifted) Macdonald polynomials: $q$-Integral
representation and combinatorial formula},
Compositio Math.\  \textbf{112}, 1998, no.\ 2, 147--182.




\bibitem{Ok3}
A.~Okounkov,
\emph{Binomial formula for Macdonald polynomials
and applications},
Math.\ Res.\ Lett.\ \textbf{4}, 1997, 533-553.



\bibitem{Ok4}
A.~Okounkov,
\emph{$BC_n$-type shifted Macdonald polynomials and binomial formula for
Koornwinder polynomials},
Transform.\ Groups \textbf{3}, 1998, no.\ 2, 181--207.

\bibitem{Ok5}
A.~Okounkov,
\emph{A characterization of interpolation
Macdonald polynomials},  Adv.\ in Appl.\ Math.\  \textbf{20}, 1998, no.\ 4,
395--428.




\bibitem{OO}
 A.~Okounkov and G.~Olshanski,
\emph{Shifted Schur functions},
 Algebra i Analiz (St.~Petersburg Math. J.)
\textbf{9} 
no.\ 2, 1997, 73--146. 




\bibitem{OO2}
 A.~Okounkov and G.~Olshanski,
\emph{ Shifted Jack polynomials, binomial formula,
and applications}, Math.\ Res.\ Lett.\ \textbf{4}, 1997, no.\ 1, 69--78.


\bibitem{OO3}
 A.~Okounkov and G.~Olshanski,
\emph{Shifted Schur functions II}
A.~A.~Kirillov Seminar on Representation Theory, 
edited by  G.~Olshanski,  
American Mathematical Society Translations, Ser.~2,
Amer.\ Math.\ Soc., Providence, 1997. 



\bibitem{OO4}
 A.~Okounkov and G.~Olshanski,
\emph{Asymptotics of Jack polynomials as 
the number of variables goes to infinity}, 
Internat.\ Math.\ Res.\ Notices 1998, no.\ 13, 641--682. 


\bibitem{S1}
 S.~Sahi,
\emph{The spectrum of certain invariant differential operators
associated to a Hermitian symmetric space}
Lie Theory and Geometry: In Honor of Bertram Kostant,
edited by J.-L.~Brylinski, R. Brylinski, V.~Guillemin, V. Kac,
Progress in Mathematics, \textbf{123}, Birkh\"auser, 1994, 569--576. 


\bibitem{S2} 
S.~Sahi,
\emph{Interpolation, integrality, and a generalization
of Macdonald's polynomials}, 
Internat.\ Math.\ Res.\ Notices 1996, no.\ 10, 457--471. 


\end{thebibliography}
\end{document}